\RequirePackage{ifpdf}
\ifpdf 
\documentclass[pdftex]{sigma}
\else
\documentclass{sigma}
\fi

\numberwithin{equation}{section}
\numberwithin{theorem}{section}
\numberwithin{proposition}{section}
\numberwithin{lemma}{section}
\numberwithin{corollary}{section}
\numberwithin{definition}{section}
\numberwithin{example}{section}
\numberwithin{remark}{section}
\numberwithin{note}{section}

\makeatletter

\def \Real {{\mathbb R}}
\def \Integer {{\mathbb Z}}

\def \Complex {{\mathbb C}}
\def \Rational {{\mathbb Q}}

\def\kk{\Bbbk}

\def\FF{{\mathbb F}}

\def \lnorm#1\rnorm {\vphantom{#1}\left\|\smash{#1}\right\|}
\def \lmod#1\rmod {\vphantom{#1}\left|\smash{#1}\right|}

\newcommand \bydef {\stackrel{\text{\rm def}}{=}}

\newcommand{\pder}[2] {\frac{\partial #1}{\partial #2}}

\@ifundefined{name}{\newcommand \name[1] {\mathop{\rm #1}\nolimits}}%
{\relax}

\renewcommand \phi {\varphi}
\renewcommand \rho {\varrho}

\let \ZZ=\Integer
\let \CC=\Complex

\def \c {\sp{(c)}}

\def \c {\sp{(c)}}
\def \Euler {\mathop{\mathcal {E}}}
\def \Laplace {L}
\def \bare {{\bar e}}

\begin{document}

\allowdisplaybreaks

\renewcommand{\thefootnote}{$\star$}

\renewcommand{\PaperNumber}{057}

\FirstPageHeading

\ShortArticleName{Dunkl Operators and Canonical Invariants of Ref\/lection Groups}

\ArticleName{Dunkl Operators and Canonical Invariants\\ of Ref\/lection Groups\footnote{This paper is a contribution to the Special
Issue on Dunkl Operators and Related Topics. The full collection
is available at
\href{http://www.emis.de/journals/SIGMA/Dunkl_operators.html}{http://www.emis.de/journals/SIGMA/Dunkl\_{}operators.html}}}

\Author{Arkady BERENSTEIN~$^\dag$ and Yurii BURMAN~$^{\ddag\S}$}

\AuthorNameForHeading{A.~Berenstein and Yu.~Burman}

\Address{$^\dag$~Department of Mathematics, University of Oregon,
Eugene, OR 97403, USA}
\EmailD{\href{mailto:arkadiy@math.uoregon.edu}{arkadiy@math.uoregon.edu}}

\Address{$^\ddag$~Independent University of Moscow, 11
B.~Vlassievsky per., 121002  Moscow, Russia}
\EmailD{\href{mailto:burman@mccme.ru}{burman@mccme.ru}}

\Address{$^\S$~Higher School of Economics, 20 Myasnitskaya Str.,  101000 Moscow, Russia}

\ArticleDates{Received December 14, 2008, in f\/inal form May 21,
2009; Published online June 03, 2009}

\Abstract{Using Dunkl operators, we introduce a continuous family of
canonical invariants of f\/inite ref\/lection groups. We verify that the
elementary canonical invariants of the symmetric group are deformations of
the elementary symmetric polynomials. We also compute the canonical
invariants for all dihedral groups as certain hypergeometric functions.}


\Classification{20F55; 15A72}


\section{Introduction and main results}

Let $V_\Real$ be a real vector space with a scalar product, and $W\subset
O(V_\Real)$ be a f\/inite group generated by ref\/lections. In this paper we construct a family of $W$-invariants (which
we refer to as {\em canonical invariants}) in $S(V)$, where $V = \Complex
\otimes V_\Real$, by means of {\em Dunkl operators} (see \cite{First}).
These canonical invariants form a basis in $S(V)^W$ (depending on a
continuous parameter $c$) and, as such, include both the {\em
$c$-elementary} and {\em $c$-quasiharmonic} invariants introduced in our
earlier paper \cite{TAMS}. Using this technique, we prove that for $W =
S_n$ the $c$-elementary invariants are deformations of the elementary
symmetric polynomials in the vicinity of $c=1/n$.

Dunkl operators $\nabla_y$, $y\in V^*$, are dif\/ferential-dif\/ference
operators f\/irst introduced by Charles Dunkl in \cite{First} and given (for
any $W$) by:
 \begin{gather*}
\nabla_y = \partial_{y} - \sum_{s\in S} c(s)\frac{\langle
y,\alpha_s\rangle}{\alpha_s} (1 - s),
 \end{gather*}
where $S$ is the set of all ref\/lections in $W$, $c:S\to \Complex$ is a
$W$-invariant function on $S$, and $\alpha_s\in V$ is the root of the
ref\/lection $s$. In particular, for $W = S_n$,
 \begin{gather*}
\nabla_y = \partial_{y} - c\sum_{1 \le i < j \le n} \frac{y_i - y_j}{x_i -
x_j}(1 - s_{ij}),
 \end{gather*}
where $s_{ij} \in S_n$ is the transposition switching $x_i$ and $x_j$.

The remarkable result by Charles Dunkl that all $\nabla_y$ commute allows
to def\/ine the opera\-tors~$\nabla_p$, $p\in S(V^*)$, by
$\nabla_{p+q}=\nabla_p+\nabla_q, \nabla_{pq}=\nabla_p\nabla_q$ for all
$p,q\in S(V^*)$.

In what follows we will mostly think of $c$ as a formal parameter in the
af\/f\/ine space ${\mathbb A}^{S/W}$, where $S/W$ is the set of $W$-orbits in
$S$. Using the notation $S_c(V) = \Complex(c)\otimes S(V)$, consider each
Dunkl operator $\nabla_p$ as a $\Complex$-linear map $S(V) \to S_c(V)$ (or,
extending scalars, as a $\Complex(c)$-linear endomorphism of $S_c(V)$).

Thus, the association $p \mapsto \nabla_p$ def\/ines an action of $S_c(V^*)$
on $S_c(V)$ by dif\/ferential-dif\/ference operators. In turn, this action and
the isomorphism $V \cong V^*$ (hence $S_c(V^*) \cong S_c(V)$) given by the
scalar product def\/ine the bilinear form $(\cdot,\cdot)_c: S_c(V) \times
S_c(V) \to \Complex(c)$ by
 \begin{gather}\label{Eq:DefForm}
(f,g)_c \bydef [\nabla_f(g)]_0,
 \end{gather}
for all $f,g \in S_c(V)$ where $x \mapsto [x]_0$ is the constant term
projection $S_c(V) \to \Complex(c)$. Clearly, $(f,g)_c = 0$ if $f$, $g$ are
homogeneous and $\deg f \ne \deg g$.

The form \eqref{Eq:DefForm} is symmetric and its specialization at generic
$c:S\to \CC$ and $c=0$ is nondegenerate. Understanding the values of $c$
when the specialization of the form is degenerate and the structure of the
radical  is crucial for the study of representations of the rational
Cherednik algebra $H_c(W)$ (see e.g.\ \cite{DuInter,DunklJeuOpdam,BEG}).

A classical Chevalley theorem \cite{Chev} says that the algebra $S(V)^W$ of
$W$-invariants in $S(V)$ is isomorphic to the algebra of polynomials
$\Complex[u_1, \dots, u_\ell]$ of certain homogeneous elements $u_1, \dots,
u_\ell$, where $\ell \bydef \dim V$. Throughout the paper we will call such
$u_1, \dots, u_\ell$ {\em homogeneous generators} or, collectively, a {\em
homogeneous generating set} of $S(V)^W$. The homogeneous generators $u_1,
\dots, u_\ell$ are not unique, but their degrees $d_1, \dots, d_\ell$
(which we traditionally list in the increasing order) are uniquely def\/ined
for each group $W$; they are called the exponents of the group. In
particular, $d_1 = 2$ if\/f $V^W = \{0\}$; the largest exponent $h \bydef
d_\ell$ is called the Coxeter number of $W$. The monomials $u^a \bydef
u_1^{a_1} \cdots u_\ell^{a_\ell}$ where $a_1, \dots, a_\ell \in
\Integer_{\ge 0}$ form an additive basis in $S(V)^W$.

Let $\prec$ be the {\em inverse lexicographic order} on $\Integer_{\ge
0}^\ell$: for $a,a'\in \Integer_{\ge 0}^\ell$ we write $a'\prec a$ if the
{\em last} non-zero coordinate of the vector $a-a'$ is positive. The
following is our f\/irst result asserting the existence and uniqueness of
canonical invariants:

 \begin{theorem}[Canonical invariants] \label{Th:CanonBasis}
Suppose that the degrees $d_1, \dots, d_\ell$ are all distinct. Then for
each $a = (a_1,\dots,a_\ell) \in \Integer_{\ge 0}^\ell$ there exists a
homogeneous element $b_a = b_a\c \in S_c(V)^W = \Complex(c) \otimes S(V)^W$
unique up to multiplication by a complex constant and such that for any
homogeneous generating set $u_1, \dots, u_\ell$ of $S(V)^W$ one has:
\renewcommand{\theenumi}{$\arabic{enumi}$}
\begin{enumerate}\itemsep=0pt
\item\label{It:Lead} $b_a \in \Complex^\times \cdot u^a + \sum_{a' \prec a}
\Complex(c)\cdot u^{a'}$;

\item\label{It:Othog} $(u^{a'},b_a)_c = 0$ whenever $a' \prec a$.
\end{enumerate}
\end{theorem}

We will prove Theorem \ref{Th:CanonBasis} in Section \ref{subsec:Basis}. We
will refer to each element $b_a$ as a {\em canonical $W$-invariant} in
$S_c(V)^W$ and to the set ${\mathbf B} = \{b_a \mid a\in \Integer_{\ge
0}^\ell\}$, as the {\em canonical basis} of $S_c(V)^W$. By the
construction, the canonical basis ${\mathbf B}$ is orthogonal with respect
to the form \eqref{Eq:DefForm}.

 \begin{remark}\label{rem:equal exponents}
We can extend the theorem to the case when $d_k=d_{k+1}$ for some $k$. If
$V$ is irreducible, then this happens only when $W$ is of type $D_\ell$
with even $\ell$ and $k = \ell/2$. In this case $V=\sum_{i=1}^\ell \Complex
\cdot x_i$, the positive roots are of the form $x_i\pm x_j$, and let
$\sigma:V\to V$ be the involution given by $\sigma(x_i)=\begin{cases} x_i
&\text{if \(i<\ell\)}\\ -x_\ell&\text{if \(i=\ell\)}\\ \end{cases}$, i.e.,
$\sigma$ is acting on roots as the symmetry of the Dynkin diagram. Then
Theorem \ref{Th:CanonBasis} holds verbatim for any choice of homogeneous
generators $u_1, \dots, u_\ell$ of $S(V)^W$ such that
$\sigma(u_{\ell/2})=-u_{\ell/2}$ (i.e., $u_{\ell/2}\in \Complex \cdot
x_1\cdots x_\ell$) and $\sigma(u_j)=u_j$ for all $j\ne \ell/2$.

An equality $d_k=d_{k+1}$ can also happen when  $V$ is reducible, i.e.,
$V=V_1\oplus V_2$, $W=W_1\times  W_2$ and each $W_i$ is a ref\/lection group
of $V_i$. This case can be handled by induction because
$S(V)^W=S(V_1)^{W_1}\otimes S(V_2)^{W_2}$.
 \end{remark}

 \begin{remark} \label{Rm:Complex}
Theorem \ref{Th:CanonBasis} generalizes to all complex ref\/lection groups if
one replaces the symmetric bilinear form on $V$ with the Hermitian one that
canonically extends the $W$-invariant Hermitian form on $S_c(V)$ (provided
that $c(s^{-1})=\overline {c(s)}$ for all complex ref\/lections $s$). The
case of equal degrees $d_k$ can be treated along the lines of Remark
\ref{rem:equal exponents}. More precisely, the phenomenon $d_k = d_{k+1}$
occurs only for the following irreducible complex ref\/lection groups (see
e.g., \cite{Cohen,BrMaRo}):
 \begin{enumerate}\itemsep=0pt
\item\label{It:Series} The series $G(m,p,\ell)$ with $\ell \ge 2$,
$p|\ell$, $p|m$, and $d_k=d_{k+1}=m\ell/p$, $k=\ell/p$.

\item\label{It:Except} The exceptional groups $G_7$, $G_{11}$, $G_{19}$ of
rank $\ell = 2$ with $d_1 = d_2 = 12, 24, 60$, respectively.
 \end{enumerate}

In the case \ref{It:Series}, similarly to Remark \ref{rem:equal exponents},
one has $V = \sum_{i=1}^\ell \Complex \cdot x_i$, $\sigma: V \to V$ is the
automorphism given by $\sigma(x_i) = \begin{cases} x_i &\text{if
\(i<\ell\)}\\ \zeta x_\ell & \text{if \(i=\ell\)}\\ \end{cases}$, where
$\zeta$ is an $m$-th primitive root of unity. Then Theorem
\ref{Th:CanonBasis} holds verbatim for any choice of homogeneous generators
$u_1, \dots, u_\ell$ of $S(V)^W$ such that $\sigma(u_k) = \zeta ^{m/p} u_k$
(i.e., $u_k\in \Complex \cdot (x_1\cdots x_\ell)^{m/p}$) and $\sigma(u_i) =
u_i$ for all $i\ne k$.

In the case \ref{It:Except} one can use various embeddings of rank $2$
complex ref\/lection groups (see e.g., \cite[Section 3]{Cohen}) to acquire
canonical invariants. For instance, $G_5$ is a normal subgroup of index $2$
in $G_7$ and $G_5$ has degrees $(6,12)$, which implies that if
$\{b_{(a_1,a_2)}\c \mid a_1,a_2\in \ZZ_{\ge 0}\}$ is the canonical basis
for $S(V)^{G_5}$, then the set $\{b_{(2a_1,a_2)}\c|a_1,a_2\in \ZZ_{\ge
0}\}$ is a (canonical) basis of $S(V)^{G_7}$.

Therefore, the invariants $b\c_a$ make sense for  all complex ref\/lection
groups.
 \end{remark}

Assume that $V^W = \{0\}$, i.e.\ $d_1=2$, and denote by $L$ the {\em Dunkl
Laplacian}  $\nabla_{e_2}:S_c(V) \to S_c(V)$, where $e_2$ is the only (up
to a scalar multiple) quadratic $W$-invariant in $S(V)$. Clearly, the
restriction of $L$ to $S_c(V)^W$ is a well-def\/ined linear operator
$S_c(V)^W\to S_c(V)^W$.

 \begin{proposition}\label{pr:kernel of sl_2}
Assume that $V^W = \{0\}$. Then
\renewcommand{\theenumi}{$(\alph{enumi})$}
\renewcommand{\labelenumi}{\theenumi}
 \begin{enumerate}\itemsep=0pt
\item\label{It:Kernel} For each $r \ge 0$ the span of all
$b_{(a_1,a_2,\dots,a_\ell)}$ with $a_1 \le r$ is the kernel of the operator
$\left. L^{r+1} \right|_{S_c(V)^W}$.

\item\label{It:MultE2} For each $a=(a_1,\dots,a_\ell)$ we have:
 \begin{gather*}
b_a = e_2^{a_1} b_{(0,a_2,\dots,a_\ell)}
 \end{gather*}
$($in particular, $b_{(a_1,0,\dots,0)} = e_2^{a_1})$.
 \end{enumerate}
 \end{proposition}

See Section \ref{subsec:elementary invariants} for the proof.

The elements $b_{(0,a_2,\dots,a_\ell)}$ of the canonical basis are more
elusive, however we compute them completely when $W$ is a dihedral group.

 \begin{theorem} \label{Th:DihCconst}
Let $W = I_2(m)$ be the dihedral group of order $2m$, $V=\Complex^2$.
 \renewcommand{\theenumi}{$(\alph{enumi})$}
\renewcommand{\labelenumi}{\theenumi}
 \begin{enumerate}\itemsep=0pt
\item\label{It:Cconst} If  $c(s_1)=c(s_2)=c$, then the generating function
of all $b_{(0,k)}$ is given by
 \begin{gather*}
\sum_{k\ge 0} \binom{c}{k} b_{(0,k)}t^k = \left(1+e_m t+e_2^m t^2\right)^c,
 \end{gather*}
where $e_2$ and $e_m$ are elementary $W$-invariants $($of degrees $2$ and $m$
respectively$)$.

\item\label{It:CNotConst} If $m$ is even and $c(s_1)\ne c(s_2)$, then
$($using the notation $C \bydef c(s_1)+c(s_2)$, $\delta \bydef
c(s_2)-c(s_1)$, $e'_m = \frac{1}{4} e_m - \frac{1}{2} e_2^{m/2})$ we have:
 \begin{gather}\label{eq:integral generating series m even}
\sum_{k\ge 0} \frac{\Gamma(\frac{C-\delta-1}{2}) \Gamma(2k-C)} {\Gamma(k -
\frac{C+\delta-1}{2})}  b_{(0,k)} t^k  = \int_0^1\!
\left(1-\tau+t\tau(e_2^{m/2} +\tau
e'_m)\right)^{\frac{C-\delta-1}{2}}\tau^{-C-1} d\tau.\!\!
 \end{gather}
 \end{enumerate}
 \end{theorem}

We will prove the theorem in Section \ref{sect:dihedral} by explicitly
reducing the Dunkl Laplacians to the Jacobi operators. In fact, it is easy
to see that the formula \eqref{eq:integral generating series m even} is
equivalent to:
 \begin{gather*}
b_{(0,k)} = \frac{k! e_2^{mk/2} }{4^k \binom{2k-C+1}{k}}
P_k^{\big(-\frac{C+\delta+1}{2},-\frac{C-\delta+1}{2}\big)}\left(\frac{e_m}{2
e_2^{m/2}}\right)
 \end{gather*}
where $P_k^{(a,b)}(y)$ is the $k$-th Jacobi polynomial (see e.g.\
\cite[Section 6.3]{SpecFunc} or formula \eqref{eq:Jacobi polynomial}
below). This and other of our arguments bear some similarity with methods
of the seminal papers \cite{First} and \cite{DunklJeuOpdam} where Jacobi
polynomials were f\/irst studied in the context of Dunkl operators.

Returning to the general case, note that $\deg b_a = \sum d_k a_k$. For each
$d \in \Integer_{\ge 0}$ such that $S(V)^W_d\ne \{0\}$ we set $e_d\c \bydef
b_{a_{\max}}$, where $a_{\max}\in \Integer_{\ge 0}^\ell$ is maximal with
respect to $\prec$ among all $a\in \Integer_{\ge 0}^\ell$ such that $\sum
d_k a_k=d$. By the construction, $\deg e_d\c = d$. The following result was
essentially proved in our previous paper \cite{TAMS}.
 \begin{theorem} \label{Th:Elem}
\renewcommand{\theenumi}{$(\alph{enumi})$}
\renewcommand{\labelenumi}{\theenumi}
Let the exponents $d_1 < \dots < d_\ell$ be pairwise distinct. Then
 \begin{enumerate}\itemsep=0pt
\item\label{It:Generates} The elements $e_{d_1}\c,\dots,e_{d_\ell}\c$
generate the algebra $S_c(V)^W$.

\item\label{It:EqForEDk} Each $e_{d_k}\c$ is determined $($up to a multiple$)$
by its homogeneity degree $d=d_k$ and the equation $\nabla_P (e_{d_k}\c) =
0$ for any $W$-invariant polynomial $P \in S_c(V)^W$ such that $\deg P <
d_k$.

\item\label{It:EqForKh} For each $k = 1, 2, \dots$ there is a unique, up to
a multiple, element $e_{kh}\c \in S_c(V)^W$ $($where $h = d_\ell$ is the
Coxeter number$)$ of the homogeneity degree $d = kh$ satisfying the equation
$\nabla_P (e_{kh}\c) = 0$ for any $W$-invariant polynomial $P \in S(V)^W$
such that $\deg P < h$.
 \end{enumerate}
 \end{theorem}

We will give a new proof of Theorem \ref{Th:Elem} in Section
\ref{subsec:elementary invariants}. The proof will rely on the construction
of canonical invariants in Theorem \ref{Th:CanonBasis}.

Following \cite{TAMS}, we refer to each $e_{d_k}\c$ as the {\em canonical
elementary $W$-invariant} and each $e_{kh}\c$ as the {\it canonical
quasiharmonic $W$-invariant}.

The elementary invariants for $c = 0$ were, most apparently, def\/ined by
Dynkin (see e.g.\ \cite{Kos}) and later explicitly computed by K.~Iwasaki
in \cite{Iwasaki}. We extend the results of \cite{Iwasaki} to all $c$ in
Theorem \ref{Th:Iwasaki} below.

We will also construct elementary invariants for $W = S_n$, $V =
\Complex^n$ with the natural $S_n$-action. It is convenient to identify
$S_c(V)$ with the algebra $\Complex(c)[x_1,\ldots,x_n]$ of polynomials in
$n$ variables depending rationally on $c$. The degrees $d_k$ are here $d_k
= k$, $k = 1, \dots, n$, so Theorem \ref{Th:CanonBasis} and Theorem
\ref{Th:Elem} are applicable.

To give the explicit formula for the invariants $e_k\c$ def\/ine polynomials
$\mu_k\c \in \Complex(c)[x_1,\ldots,x_n]$, $k = 2, \dots, n$, by
 \begin{gather*}
\mu_k\c = \sum_{s=1}^k (-1)^s x_s (\Delta(\nabla_{x_1}, \dots,
\widehat{\nabla_{x_s}}, \dots, \nabla_{x_k})) \Delta(x_1,\ldots,x_k),
 \end{gather*}
where $\Delta(z_1,\dots,z_r) = \prod_{1 \le i<j \le r} (z_i-z_j)$ is the
Vandermonde determinant. Clearly, $\mu_k\c \in
\Complex(c)[x_1,\dots,x_n]^{S_k\times S_{n-k}}$.

 \begin{theorem} \label{Th:Iwasaki}
For all $2 \le k \le n$ there exists  $\alpha_{k,n} \in \Complex(c)^\times$
such that the $k$-th elementary canonical invariant $e_k\c \in
\Complex(c)[x_1, \dots, x_n]^{S_n}$ is given by:
 \begin{gather}\label{Eq:Iwa}
e_k\c = \alpha_{k,n}(c) \sum_{w\in S_n/(S_k\times S_{n-k})}
w\big(\mu_k\c\big).
 \end{gather}
 \end{theorem}

We prove Theorem \ref{Th:Iwasaki} in Section \ref{subsec:Iwasaki}. Our
proof (as well as the formula \eqref{Eq:Iwa}) is very similar to the one by
K.~Iwasaki who (using $\partial_p$ instead of Dunkl operators $\nabla_p$)
computed $e_k^{(0)}$ in \cite{Iwasaki}. Following his argument, one can
construct the elementary canonical invariants $e_{d_k}\c$ for other
classical groups as well.

Note that the formula \eqref{Eq:Iwa} resembles the polynomial expansion of
the elementary symmetric polynomial $e_k = e_k(x_1,\dots, x_n)$:
 \begin{gather*}
e_k = \sum_{w\in S_n/(S_k\times S_{n-k})}
w(x_1\cdots x_k)=\sum_{1 \le j_1 < \dots < j_k \le n} x_{j_1} \cdots x_{j_k}.
 \end{gather*}
The following main result demonstrates that this observation is not a mere
coincidence.

 \begin{theorem} \label{Th:Main}
Let $W = S_n$. Then for all $k = 2, \dots, n$ the elementary canonical
invariants $e_k\c$ have no poles at the singular value $c = 1/n$, and
 \begin{gather*}
\lim_{c\to 1/n} e_k\c = e_k\left(x_1-\frac{e_1(x)}{n}, \dots,
x_n-\frac{e_1(x)}{n}\right).
 \end{gather*}
 \end{theorem}

This result allows to introduce the elementary invariant polynomials for
other ref\/lection groups via $e_{d_k} = \lim\limits_{c\to 1/h} e_{d_k}\c$, where
$h$ is the Coxeter number.

We will prove Theorem \ref{Th:Main} in Section \ref{SSec:PrMain} by
analyzing the behaviour of the form \eqref{Eq:DefForm} near $c = 1/n$.
Note, however, that we could not derive the theorem directly from the
explicit formula~\eqref{Eq:Iwa}.

 \begin{example}
Denote $\bare_k(x) \bydef e_k\bigl(x_1-\frac{e_1(x)}{n}, \dots,
x_n-\frac{e_1(x)}{n}\bigr)$. It is easy to see that $e_1\c = 0$, $e_k\c =
\bare_k$ for $k=2,3$. Direct computations for all $n$ using \eqref{Eq:Iwa}
show that
 \begin{gather}
e_4\c = \frac{(n-2)(n-3)}{2n}\frac{1-nc}{(n^2-n)c - n - 1} \bare_2^2 +
\bare_4,\label{Eq:e4}\\
e_5\c = \frac{(n-3)(n-4)}{n} \frac{1-nc}{(n^2-n)c - n - 5} \bare_2 \bare_3
+ \bare_5,\label{Eq:e5}
 \end{gather}
thus conf\/irming Theorem \ref{Th:Main}.
 \end{example}

 \begin{remark}
The def\/inition of the canonical invariants $b_a$ and some later formulas
involving them (e.g.\ \eqref{Eq:Frob} and \eqref{Eq:DefPhi}) suggest, for
$W = S_n$, a close relation between canonical invariants~$b_a$ and Jack
polynomials $J_\lambda^{(\alpha)}$ (see e.g.\ \cite{Opd} for def\/inition).
Direct computations show, though, that these polynomials are {\em not} the
same. \cite[equation (7)]{Lap} shows, in particular, that the expression of~$J_\lambda^{(\alpha)}$ via elementary symmetric polynomials $e_i$ does not
depend on $n$; for instance, $J_{(11\dots1)}^{(\alpha)} = e_k$ for all $n$
and $k$ (the partition contains $k$ units). Formulas for $b_a$, on the
contrary, contain $n$ explicitly (see e.g.~\eqref{Eq:e4}). So, the
relation between $b_a$ and Jack polynomials is yet to be clarif\/ied.
 \end{remark}

\section{The Dunkl Laplacian and the scalar product}

Throughout the section we assume that $2 = d_1 < \dots < d_\ell = h$ and
denote
 \begin{gather*}
e_2=\sum_{i=1}^\ell x_i^2,
 \end{gather*}
where $x_1, \dots, x_\ell$ is any orthonormal basis in the real space
$V_\Real$. Obviously, $e_2$ is a unique (up to a scalar multiple) quadratic
$W$-invariant in $S^2(V)$. The operator $\Laplace = \nabla_{e_2} = \sum_i
\nabla_{x_i}^2$, called the {\em Dunkl Laplacian}, is independent of the
choice of the basis $x_i$; it equals the ordinary Laplacian if $c = 0$.

The operator $L$ plays a key role in the theory of Dunkl operators for $W$.
As the following result shows, an action of any Dunkl operator can be
expressed via $L$:

 \begin{lemma}[\protect{\cite[equation (1.9)]{Berest}}]\label{le:Berest}
For any $p \in S^d(V)$ one has
 \begin{gather*}
\nabla_p = \frac{1}{d!}(\name{ad} L)^d(p) = \sum_{k=0}^d
\frac{(-1)^k}{k!(d-k)!} L^{d-k}\cdot p \cdot L^k,
 \end{gather*}
where $p$ in the right-hand side means the operator of multiplication by
$p$.
 \end{lemma}

Denote by $\Euler:S_c(V) \to S_c(V)$ the Euler vector f\/ield given by
$\Euler(f) = Nf$ for any $f \in S_c^N(V)$. Also denote $h_c \bydef
\frac{2}{\ell}\sum_{s\in S} c(s)$ (in particular, if all $c(s)$ are equal
to a single $c$, then $h_c = hc$).

 \begin{proposition}[\cite{Hec}]\label{Pp:SL2}
The operator $E$ of multiplication by $e_2$, Dunkl Laplacian $\Laplace$,
and the operator $H \bydef 2\ell(1-h_c) +4\Euler$ form a representation of
$\mathfrak{sl}_2$, that is,
 \begin{gather*}
[E,L]=H, \qquad [H,E]=2E, \qquad [H,L]=-2L.
 \end{gather*}
In particular,
 \begin{gather*}
[\Laplace, E^k] = 4k E^{k-1} (\ell(1-h_c)/2 + k - 1 + \Euler)
 \end{gather*}
for all $k\ge 0$.
 \end{proposition}

Denote $U_d=U_d\c \bydef \name{Ker} L\cap S_c^d(V) = \{f \in S^d_c(V) \mid
\Laplace(f) = 0\}$.

 \begin{lemma}\label{le:Orthog}
One has
 \begin{gather}\label{Eq:Decomp}
S_c(V) = \bigoplus_{k,d \in \Integer_{\ge 0}} e_2^k \cdot U_d\c,
 \end{gather}
where the direct summands are orthogonal with respect to $(\cdot,
\cdot)_c$. In particular, the restriction of $(\cdot, \cdot)_c$ to each
$e_2^k \cdot U_d$ is nondegenerate.
 \end{lemma}

 \begin{proof}
First, note that $S_c(V)$ is an $\mathfrak{sl}_2$-module, locally f\/inite
with respect to $L$, and $\oplus_{d\ge 0} U_d\c$ is the highest weight
space, so that decomposition \eqref{Eq:Decomp} takes place.

Furthermore, note that the operator $H$ from Proposition \ref{Pp:SL2} is
scalar on the space of polynomials of any given degree and therefore
self-adjoint; the operators $L$ and $E$ are adjoint to one another with
respect to $(\cdot, \cdot)_c$. Therefore, for $k_1\le k_2$, $d_1,d_2\ge 0$
one has
 \begin{gather*}
\big(e_2^{k_1}\cdot U_{d_1}\c,e_2^{k_2}\cdot U_{d_2}\c\big)_c  =
\big(E^{k_1}\big(U_{d_1}\c\big), E^{k_2}\big(U_{d_2}\c\big)\big)_c = \big(L^{k_2}E^{k_1}\big(U_{d_1}\c\big),
U_{d_2}\c\big)_c\\
\phantom{\big(e_2^{k_1}\cdot U_{d_1}\c,e_2^{k_2}\cdot U_{d_2}\c\big)_c}{} =\delta_{k_1,k_2}\cdot \big(U_{d_1}\c, U_{d_2}\c\big)_c = \delta_{k_1,k_2}
\delta_{d_1,d_2} \cdot \Complex(c).
 \end{gather*}
This proves the orthogonality of the decomposition. In particular, this
implies that the restriction of the nondegenerate form $(\cdot, \cdot)_c$
to each  $e_2^k \cdot U_d$ is nondegenerate. The lemma is proved.
 \end{proof}

Using this, we compute the form $(\cdot,\cdot)_c$ as follows. Denote by
$\phi_c$ a (unique) linear function $S_c(V)\to \Complex(c)$ such that:{\samepage
 \begin{itemize}\itemsep=0pt
\item $\phi_c(fe_2) = \phi_c(f)$ for all $f\in S_c(V)$;

\item $\phi_c(f)=[f]_0$ for all $f\in \name{Ker} L$, where $[\,\cdot\,]_0:
S_c(V) \to \CC(c)$ is the projection def\/ined in~\eqref{Eq:DefForm}.
 \end{itemize}}

 \begin{proposition}\label{Pp:Frob}
We have:
\renewcommand{\theenumi}{$(\alph{enumi})$}
\renewcommand{\labelenumi}{\theenumi}
\begin{enumerate}\itemsep=0pt
\item For $f \in e_2^kU_d\c$, $g \in S_c^{d+2k}(V)$  one has
 \begin{gather}\label{Eq:Frob}
(f,g)_c =  \phi_c(fg)\cdot 4^{d+k} k!\prod_{r=0}^{d+k-1} (\ell(1-h_c)/2 +
r).
 \end{gather}

\item If $c: S/W \to \Real_{< 1/2}$, then the restriction of
$\phi_c$ to $S(V_\Real)$ is given by:
 \begin{gather}\label{Eq:DefPhi}
\phi_c(f) =\frac{\int_{\Omega^{\ell-1}} f(x) \cdot \prod_{s \in S}
\lmod\alpha_s(x)\rmod^{-2c(s)}\,dx} {\int_{\Omega^{\ell-1}} \prod_{s \in S}
\lmod\alpha_s(x)\rmod^{-2c(s)}\,dx},
 \end{gather}
where $\alpha_s\in V_\Real$ is a coroot of a reflection $s \in S$,
$\Omega^{\ell-1} = \{x \in V_\Real \mid e_2(x) = 1\}$ is the unit sphere
in $V_\Real$, and an element $f\in S(V)$ is identified with a polynomial on
$V$.
 \end{enumerate}
 \end{proposition}

 \begin{proof}
Assume f\/irst that the function $c$ takes only negative real values and
def\/ine $\phi_c$ by equation~\eqref{Eq:DefPhi}. Now if $k = 0$ and $f,g\in
U_d\c$, then the result follows from \cite[Theorem 5.2.4]{DXu}.

By definition, $\phi_c(e_2 f) = \phi_c(f)$ for any $f$. Now if $f = e_2^k
\tilde f$, $g = e_2^{k'} \tilde g$ where $\tilde f \in U_d\c$, $\tilde g
\in U_{d'}\c$ with $d \ne d'$ then $(\tilde f, \tilde g)_c = 0$, hence
$\phi_c(\tilde f \tilde g) = 0$ and therefore $\phi_c(fg) = 0$. So taking
$f \in e_2^kU_d\c$ and $g = \sum_{r}e_2^r \tilde g_r \in S_c^{d+2k}(V)$,
where $\tilde g_r\in U_{d+2k-2r}$, we see that
 \begin{gather*}
\phi_c(fg)=\phi_c(\tilde f\cdot \tilde g_k).
 \end{gather*}

On the other hand, decomposition \eqref{Eq:Decomp} guarantees that $(f,g)_c
= (\tilde f,\tilde g_k)_c$. Therefore, to verify \eqref{Eq:Frob} for any
$k$ it suf\/f\/ices to take $g = e_2^k\tilde g$ for $g \in U_d\c$.

Assume that $k>0$. Then Proposition \ref{Pp:SL2} implies that
 \begin{gather*}
(e_2^k\tilde f, e_2^k\tilde g)_c  = (E^k(\tilde f),E^k(\tilde g))_c =
(E^{k-1}(\tilde f), LE^k(\tilde g))_c = (E^{k-1}(\tilde f),[L,E^k](\tilde
g))_c\\
\phantom{(e_2^k\tilde f, e_2^k\tilde g)_c}{} =(E^{k-1}(\tilde f), 4k E^{k-1} (\ell(1-h_c)/2 + k - 1 +
\Euler)\tilde g))_c\\
\phantom{(e_2^k\tilde f, e_2^k\tilde g)_c}{}=4k(\ell(1-h_c)/2 + k - 1 +d)(e_2^{k-1}(\tilde f),e_2^{k-1}\tilde
g)_c.
 \end{gather*}
Therefore, by induction on $k$,
 \begin{gather*}
(e_2^k\tilde f,e_2^k\tilde g)_c = (\tilde f,\tilde g)_c\cdot
\prod_{r=1}^k 4r(\ell(1-h_c)/2 + r - 1 +d),
 \end{gather*}
which f\/inishes the proof for $c$ negative real. Now \eqref{Eq:Frob} implies
that for $c$ negative real the value $(f,g)_c$ depends only on the product
$fg$ (provided $d$ and $k$ are f\/ixed). Since $(f,g)_c$ is a rational
function of the values of $c$, this holds true for all $c$ as well -- so,
one can use \eqref{Eq:Frob} to def\/ine $\phi_c$ in the general case.
 \end{proof}

The following is the main result of the section. Denote by $S(V)_+$ the
kernel of the constant term projection $u\to [u]_0$, see
\eqref{Eq:DefForm}. Def\/ine a symmetric bilinear form $\Phi_c: S(V)_+ \times
S(V)_+ \to \Complex(c)$ as $\Phi_c(u,v) \bydef (u,v)_c/(1-h_c)$. This form
extends naturally to $\Complex[c] \otimes S(V)$. Def\/ine now the form
$\overline \Phi_c$ on $\Complex[c]/(1-h_c) \otimes S(V)$ taking values in
$\Complex[c]/(1-h_c)$ by
 \begin{gather*}
\overline \Phi_c(u,v) = \pi(\Phi_c(\tilde u,\tilde v)),
 \end{gather*}
where $\pi: \Complex[c] \to \Complex[c]/(1-h_c)$ is the canonical
projection and $\tilde u, \tilde v \in S_c(V)$ are any elements such that
$u = \pi(\tilde u)$, $v = \pi(\tilde v)$.

 \begin{theorem}\label{Th:NonDeg}\qquad
\renewcommand{\theenumi}{$(\alph{enumi})$}
\renewcommand{\labelenumi}{\theenumi}
 \begin{enumerate}\itemsep=0pt
\item\label{It:NoPoles} The form $\Phi_c$ takes its values in
$\Complex[c]$.

\item\label{It:PhiNonDeg} $\overline \Phi_c(u,u)\ne 0$ for any non-zero
element $u\in \Real[c]/(1-h_c)\otimes_\Real S(V_\Real)_+$.

\item\label{It:PhiBarSub} For any $U_\Real \subset S(V_\Real)_+$ the
restriction of $\overline \Phi_c$ to $\Complex[c]/(1-h_c)\otimes  U_\Real$
is nondegenerate.
 \end{enumerate}
 \end{theorem}

 \begin{proof}
Prove \ref{It:NoPoles} by induction on the degree of
$u$. Indeed, it follows from \cite[Proposition 2.1]{BEG} that for any
$x,y\in V$ one has:
 \begin{gather*}
(x,y)_c = (1-h_c)(x,y)_0,
 \end{gather*}
where $(x,y)_0$ is the (complexif\/ied) $W$-invariant form on $V$. Therefore,
$\left.\Phi_c\right|_{V\times V} = (\cdot,\cdot)_0$. Furthermore, assume
that $\Phi_c(u,v) \in \Complex[c]$ for all $u,v\in S^{<d}(V)_+$. Then for
any $u_1\in S^{d_1}(V)$, $u_2\in S^{d_2}(V)$,  $v\in S^d(V)$, where
$d_1+d_2=d$, we have
 \begin{gather*}
\Phi_c(u_1u_2,v) = \Phi_c(u_2,\nabla_{u_1}(v)) \in \Phi_c(u_2,\Complex[c]
\otimes S^{d_2}(V)) \subset \Complex[c].
 \end{gather*}
This proves \ref{It:NoPoles}.

It is possible to prove \ref{It:PhiNonDeg} now. Let $H_0 \bydef \{c
\in \Complex^{S/W} \mid h_c = 1\}$; it is an af\/f\/ine hyperplane in the
af\/f\/ine space ${\mathbb A}^{S/W}$. Then $\Complex[H_0] =
\Complex[c]/(1-h_c)$ and $\Real[H_0] = \Real[c]/(1-h_c)$ are integral
domains (if $\lmod S/W\rmod=1$ then $H_0$ is a point $c=1/h$ and
$\Real[H_0]=\Real$).

Let ${\mathcal A} \subset \Real({\mathbb A}^{S/W})$ be the   algebra of all
real-valued rational functions on the af\/f\/ine space ${\mathbb A}^{S/W}$
regular at $H_0$. This algebra is local with the maximal ideal ${\mathfrak
m}=(1-h_c)$, and ${\mathcal A}/{\mathfrak m}$ is isomorphic to
$\Real(H_0)$, the f\/ield of fractions of $H_0$. Finally, denote $S_{\mathcal
A}(V_\Real)_+ = {\mathcal A} \otimes_\Real S(V)_+$.

 \begin{proposition} \label{pr:non-isotropic Phi}
The naturally extended ${\mathcal A}$-linear form  $\Phi_c:S_{\mathcal
A}(V_\Real)_+\times S_{\mathcal A}(V_\Real)_+\to {\mathcal A}$ satisfies:
 \begin{gather}\label{eq:non-isotropic Phi}
\text{if} \ \ \Phi_c(\tilde u,\tilde u)\in (1-h_c){\mathcal A} \ \ \text{for some} \
 \ \tilde u \in S_{\mathcal A}(V_\Real) \ \ \text{then} \ \ \tilde u\in
(1-h_c) S_{\mathcal A}(V_\Real).
 \end{gather}
 \end{proposition}

 \begin{proof}
Clearly, the $\mathfrak{sl}_2$-action from Proposition \ref{Pp:SL2}
preserves both $\Real[c]\otimes S(V_\Real)$ and $S_{\mathcal A}(V_\Real)$,
so that the orthogonal decomposition \eqref{Eq:Decomp} is valid for
$S_{\mathcal A}(V_\Real)\subset S_c(V)$.  Therefore, it suf\/f\/ices to verify
\eqref{eq:non-isotropic Phi} only for $\tilde u \in e_2^k\tilde U_d$, where
 \begin{gather*}
\tilde U_d = U_d\c\cap S_{\mathcal A}(V_\Real) = \{\tilde v \in S_{\mathcal
A}(V_\Real) \mid L(\tilde v)=0\}.
 \end{gather*}
For every such $\tilde u$ it follows from \eqref{Eq:Frob} that
 \begin{gather}\label{eq:Phi of the square}
\Phi_c(\tilde u,\tilde u) = \phi_c(\tilde u^2) \cdot 2 \cdot 4^{d+k-1}
k!\prod_{r=1}^{d+k-1} (\ell(1-h_c)/2 + r).
 \end{gather}
Since the product in the right-hand side is not divisible by $(1-h_c)$, we
see that $\phi_c(\tilde u^2)\in {\mathcal A}$ for all $\tilde u\in
S_{\mathcal A}(V_\Real)$. Implication \eqref{eq:non-isotropic Phi} is now
equivalent to the following one:
 \begin{gather}\label{eq:non-isotropic frobenius}
\text{if} \ \ \phi_c(\tilde u^2) \in (1-h_c){\mathcal A} \ \ \text{for some} \ \
\tilde u\in S_{\mathcal A}(V_\Real) \ \ \text{then} \ \ \tilde u \in
(1-h_c)S_{\mathcal A}(V_\Real).
 \end{gather}

In \eqref{eq:non-isotropic frobenius}, if $h=2$, i.e., $W = S_2$, we have
nothing to prove. Assume that $h > 2$ and let $\tilde H_0$ be the set of
all $c_0: S/W \to \Real_{<1/2}$ such that $1-h_{c_0} = 0$ and $\tilde u$
has no poles at $c_0$. By the very design, $\tilde H_0$ is a non-empty open
subset of the real hyperplane $H_0(\Real)\subset \Real^{S/W}$. Indeed, we
can write $h_c = h_1 c_1 + \dots + h_k c_k$, where $k = \lmod S/W\rmod$ and
$c_1, \ldots,c_k$ are standard coordinates on ${\mathbb A}^{S/W} = {\mathbb
A}^k$, all $h_i \in \Rational_{>0}$ and $h_1 + \dots + h_k = h$. The
intersection $H_0' = H_0(\Real) \cap (\Real_{<1/2})^k$ contains the point
$c = (1/h, \dots, 1/h)$, hence $H_0'$ is non-empty and open in
$H_0(\Real)$. The set $\tilde H_0$ is obtained from $H'_0$ by removing
poles of $\tilde u \in {\mathcal A} \otimes S(V_\Real)$; so, $\tilde H_0
\subset H_0(\Real)$ is open and non-empty as well. On the other hand, for
each $c_0 \in \tilde H_0$ the specialization $\tilde u_{c_0}\in S(V_\Real)$
of $\tilde u$ is well-def\/ined. Therefore, our choice of $\tilde u$ and
$c_0$ implies that
 \begin{gather*}
\phi_{c_0}\big(\tilde u_{c_0}^2\big) = 0.
 \end{gather*}
But the integral presentation \eqref{Eq:DefPhi} of $\phi_{c_0}$ guarantees
that $\phi_{c_0}(f^2)>0$ for all nonzero polynomials $f\in S(V_\Real)$.
Hence, $\tilde u_{c_0} = 0$, for all $c_0 \in \tilde H_0$. If $\lmod
S/W\rmod = 1$, i.e., $H_0$ is a single point $c=1/h$, then, clearly,
$c_0 = 1/h$ and $\tilde u_{1/h} = 0$ implies that $\tilde u \in (1-hc)
S_{\mathcal A}(V_\Real)$. Otherwise, if $H_0$ is at least an af\/f\/ine line,
the set $\tilde H_0$ is inf\/inite and is a set of regular points for a
rational function $c_0 \mapsto \tilde u_{c_0}$. Therefore, $\tilde u_{c_0}
= 0$ for all $c_0 \in H_0(\Real)$ and hence $\tilde u\in (1-h_c)S_{\mathcal
A}(V_\Real)$ as well, proving implications  \eqref{eq:non-isotropic
frobenius} and \eqref{eq:non-isotropic Phi}. The proposition is proved.
 \end{proof}

Part \ref{It:PhiNonDeg} of the theorem immediately follows from
Proposition \ref{pr:non-isotropic Phi}.

To prove \ref{It:PhiBarSub} we need the following obvious result:

 \begin{lemma} \label{Lm:Ext}
Let $\overline \Phi:U_0\times U_0\to \kk$ be a non-degenerate symmetric
bilinear form on a $\kk$-vector space $U_0$. Then for any field  $\FF$
containing  $\kk$ the extension of $\Phi$ to $U=\FF\otimes_\kk U_0$ is  a
non-degenerate $\FF$-bilinear form $U\times U\to \FF$.
 \end{lemma}

We will use the lemma with $\kk$ being the f\/ield of fractions of the
integral domain $\Real[c]/(1-h_c)$, $\FF$ -- the f\/ield of fractions of
$\Complex[c]/(1-h_c)$, $U_0 = \kk \otimes U_\Real$, and $\Phi=\overline
\Phi_c$. It follows from part~\ref{It:PhiNonDeg} of the theorem that the
restriction of $\overline \Phi_c$ to $U_0$ is non-degenerate. Hence Lemma~\ref{Lm:Ext} guarantees the same for $U = \FF \otimes U_0$. This proves
part~\ref{It:PhiBarSub} of the theorem.
 \end{proof}

 \begin{remark}
The proof of Theorem \ref{Th:NonDeg}\ref{It:PhiNonDeg} also implies
unitarity of $S(V_\Real)_+$ as a module over the rational Cherednik algebra
$H_{c_0}(W)$ for all $c_0\in \Real_{\le 0}^{S/W}$ and for small $c_0\in
\Real_{>0}^{S/W}$. This agrees with the results of the recent paper
\cite{ESG}, where the unitary representations of $H_c(W)$ were studied.
 \end{remark}

\section{Canonical basis and proofs of main results}

\subsection{Proof of Theorem \ref{Th:CanonBasis}}
\label{subsec:Basis}
Fix a homogeneous generating set $\{u_1, \dots, u_\ell\}$ of $S(V)^W$ and
take the basis $u^a$, $a \in \Integer_{\ge 0}^\ell$, in~$S(V)$ with the
inverse lexicographic order. Def\/ine the $\Complex$-subspaces $S_c(V)_{\prec
a}^W$ and  $S_c(V)_{\preceq a}^W$ of $S_c(V)^W$ by
 \begin{gather}\label{eq:SVpreca}
S_c(V)_{\prec a}^W = \sum_{a'\prec a}\Complex(c) u^{a'},\qquad S_c(V)_{\preceq a}^W = \Complex u^a +  S(V)_{\prec a}^W.
 \end{gather}
(clearly,  $S_c(V)_{\prec a}^W\subset S_c(V)_{\preceq a}^W$). Note f\/irst
that for each $a\in \Integer_{\ge 0}^\ell$ the spaces $S_c(V)_{\prec a}^W$,
$S_c(V)_{\preceq a}^W$ do not depend on the choice of generators $u_1,
\dots, u_\ell$ of $S(V)^W$. Indeed, let $u'_1, \dots, u'_\ell$ be another
set of generators of $S(V)^W$. Since $d_1 < d_2 < \dots < d_\ell$, one has
$u'_i = \alpha_i u_i + P_i(u_1, \dots, u_{i-1})$, where $\alpha_i \in
\Complex \setminus \{0\}$ and $P_i$ is a polynomial of $i-1$ variables for
$i = 1,2,\dots,\ell$.

We are going to def\/ine the canonical invariant $b_a\in S_c(V)_{\preceq
a}^W$ as the unique (up to a multiple) vector orthogonal to the subspace
$S_c(V)_{\prec a}^W$. However, the uniqueness of such an element requires
more arguments.

 \begin{proposition} \label{pr:non-degenerate form}
Let $U$ be any subspace of $S(V_\Real)$. Then:
\renewcommand{\theenumi}{$(\alph{enumi})$}
\renewcommand{\labelenumi}{\theenumi}
 \begin{enumerate}\itemsep=0pt
\item\label{It:Restr} the restriction of the form \eqref{Eq:DefForm} to
$U_c = \Complex(c)\otimes U$ is a non-degenerate symmetric bilinear form on
$U_c$;

\item\label{It:Unique} for any vector $u\in S(V_\Real)\setminus U$ there is
a unique $($up to a complex multiple$)$ element $b \in \Complex \cdot u+ U_c$
such that $(b,U_c)_c = 0$.
 \end{enumerate}
 \end{proposition}

 \begin{remark}
Statement \ref{It:Restr} of this proposition is essentially a Cherednik
algebra version of the statements proved in \cite[statements 5.1.20 and
5.1.21]{Mac} for double af\/f\/ine Hecke algebras.
 \end{remark}

 \begin{proof}
We need the following general fact. Let ${\mathcal A}$ be a unital
commutative local ring with no zero-divisors, ${\mathfrak m}$ its maximal
ideal, $\kk = {\mathcal A}/{\mathfrak m}$ the residue f\/ield. In what
follows we assume  that $\kk\subset \mathcal A$ so that the restriction of
the canonical projection ${\mathcal A}\to \kk$ to $\kk$ is the identity
homomorphism $\kk\to \kk$. Let $U$ be a vector space over $\kk$ and let
$\Phi:U\times U\to {\mathcal A}$ be a   $\kk$-bilinear symmetric form on
$U$; denote by $\Phi_0: U \times U \to \kk$ the residual form given by
$\Phi_0 \bydef \pi \circ \Phi$, where $\pi:{\mathcal A}\to \kk$ is the
canonical quotient map.

 \begin{lemma} \label{le:non-isotropic form}
In the notation as above assume that:
\renewcommand{\theenumi}{$\arabic{enumi}$}
 \begin{enumerate}\itemsep=0pt
\item\label{It:NonDeg} $\Phi_0(u,u) \ne 0$ for all $u\in U\setminus \{0\}$.

\item\label{It:Ideals} There exists an increasing sequence of prime ideals
 \begin{gather*}
\{0\} = {\mathfrak m}_0 \subset {\mathfrak m}_1\subset {\mathfrak m}_2
\subset \dots \subset {\mathfrak m}_k = {\mathfrak m}
 \end{gather*}
in ${\mathcal A}$ such that  ${\mathfrak m}_{i+1}/{\mathfrak m}_i$ is a
principal ideal of ${\mathcal A}/{\mathfrak m}_i$ for $i=0,1,\dots,k-1$.
 \end{enumerate}
 
Then the natural ${\mathcal A}$-bilinear extension  of $\Phi$ to
$U_{\mathcal A} \times U_{\mathcal A} \to {\mathcal A}$, where $U_{\mathcal
A}\bydef {\mathcal A}\otimes_\kk U$, satisfies $\Phi(\tilde u,\tilde u)\ne
0$ for all $\tilde u\in U_\FF\setminus \{0\}$.
 \end{lemma}

 \begin{proof}
We proceed by induction on $k$. For $k=0$, ${\bf m}=\{0\}$, $\FF =
{\mathcal A} = \kk$, and we have nothing to prove. Assume that $k \ge 1$.
Def\/ine the quotient ring ${\mathcal A}' \bydef {\mathcal A}/{\mathfrak
m}_1$, and ${\mathfrak m}'_i \bydef {\mathfrak m}_{i+1}/{\mathfrak m}_1$ in
${\mathcal A}'$ for $i = 0, \dots, k-1$. Clearly, the ring ${\mathcal A}'$
and its ideals ${\mathfrak m}'_i$ satisfy the assumptions of the lemma for
$k-1$; therefore, the inductive hypothesis holds in the following form:
 \begin{gather}\label{inductive hypothesis non-isotropic}
\text{if} \ \ \Phi(\tilde u,\tilde u) \in {\mathfrak m}_1 \ \ \text{for some} \ \
\tilde u\in {\mathcal A}\otimes U, \ \ \text{then} \ \ \tilde u\in {\mathfrak
m}_1\otimes U.
 \end{gather}

Since the ideal ${\mathfrak m}_1$ is principal, i.e., ${\mathfrak
m}_1=c_1{\mathcal A}$, we can write each non-zero vector $\tilde u\in
{\mathfrak m}_1\otimes U$ in the form $\tilde u=c_1^\ell \tilde u_0$, where
$\tilde u_0\notin {\mathfrak m}_1\otimes U$. Therefore, the equation
$\Phi(\tilde u,\tilde u)=0$ is equivalent to $\Phi(\tilde u_0,\tilde
u_0)=0$. However, applying the inductive hypothesis \eqref{inductive
hypothesis non-isotropic} to any $\tilde u_0\notin {\mathfrak m}_1\otimes
U$ satisfying $\Phi(\tilde u_0,\tilde u_0)=0$, we obtain a contradiction.
Therefore,  $\Phi(\tilde u_0,\tilde u_0)\ne 0$ for all $u_0\notin
{\mathfrak m}_1\otimes U$. Hence $\Phi(\tilde u,\tilde u)=0$ for $\tilde
u\in {\mathcal A}\otimes U$ if and only if $\tilde u=0$.

The lemma is proved.
 \end{proof}

We apply the lemma in the case when $\FF = \Complex(c) =
\Complex(c_1,\ldots,c_k)$ is the f\/ield of rational functions in the
variables $c_1,\ldots,c_k$ (where $k = \lmod S/W\rmod$ is the number of
conjugacy classes of ref\/lections in $W$), ${\mathcal A} \subset \FF$ is the
local ring of all rational functions regular at $c = 0$, and ${\mathfrak
m}_i$ is the ideal of ${\mathcal A}$ generated by $c_1, \dots, c_i$ for $i
= 0, 1, \dots, k$. Clearly, the ideals ${\mathfrak m}_i$ satisfy condition~\ref{It:Ideals} of Lemma~\ref{le:non-isotropic form}. Take $U$ to be any
subspace of $S(V_\Real)$ and let $\Phi: U \times U \to {\mathcal A} \subset
\Real(c)$ be the restriction of the form \eqref{Eq:DefForm} to $U$. Since
the specialization $\Phi_0$ of $\Phi$ at $c = 0$ is a positive def\/inite
form on $U$, condition~\ref{It:NonDeg} of Lemma \ref{le:non-isotropic form}
holds as well.

Thus, Lemma \ref{le:non-isotropic form} guarantees that for any $U\subset
S(V_\Real)$ each $\tilde u\in \Real(c)\otimes U\setminus \{0\}$ satisf\/ies
$(\tilde u,\tilde u)_c\ne 0$.

Therefore, the restriction of  the form \eqref{Eq:DefForm} to $\Real(c)
\otimes U$ is non-degenerate. By extending the coef\/f\/icients from $\Real(c)$
to $\Complex(c)$ this immediately proves assertion~\ref{It:Restr} of
Proposition \ref{pr:non-degenerate form}.

To prove assertion~\ref{It:Unique} denote $U_c^\perp = \{\tilde u'\in
\Complex(c) u + U_c \mid (\tilde u',U_c)_c = 0\}$. Clearly, $U_c^\perp \ne
0$ and $(U_c^\perp \cap U_c,U_c)_c = 0$; hence $U_c^\perp\cap U_c=0$ by
assertion~\ref{It:Restr}. This implies that $\dim U_c^\perp = 1$.
Therefore, $U_c^\perp=\Complex(c)\cdot b$ for some $b\in u + U_c$. This
completes the proof of Proposition \ref{pr:non-degenerate form}.
 \end{proof}

Now we are ready to f\/inish the proof of Theorem \ref{Th:CanonBasis}. For
each $a \in \ZZ_{\ge 0}^\ell$ denote by $S_c(V_\Real)_{\prec a} =
S_c(V_\Real) \cap S_c(V)_{\prec a}$ (see \eqref{eq:SVpreca}) the real forms
of $S(V)_{\prec a}$. Fix $u_1, \dots, u_\ell$  to be a homogeneous
generating set of $S(V_\Real)^W$ so that \eqref{eq:SVpreca} implies that
$S_c(V_\Real)_{\prec a}^W = \sum_{a'\prec a} \Real(c) u^{a'}$ and $\Complex
\otimes S_c(V_\Real)_{\prec a} = S_c(V)_{\prec a}$.

Therefore, Proposition \ref{pr:non-degenerate form} is
applicable to this situation with $U = S_c(V_\Real)_{\prec a}^W$, $u = u^a$,
and there exists a unique (up to a complex multiple) element $b = b_a\in
S_c(V)_{\preceq a}$ such that $(b_a,S_c(V)_{\prec a}^W) = 0$ (in
particular, $b_a\notin S_c(V)_{\prec a}^W$). In other words, $b_a$
satisf\/ies both conditions of Theorem \ref{Th:CanonBasis}, and the theorem
is proved.

\subsection{Proof of Proposition \ref{pr:kernel of sl_2} and Theorem
\ref{Th:Elem}}\label{subsec:elementary invariants}

 \begin{proof}[Proof of Proposition~\ref{pr:kernel of sl_2}]
To prove \ref{It:Kernel}, let $u_1(=e_2), u_2, \dots, u_\ell$ be any
homogeneous genera\-ting set of $S(V)^W$. Theorem \ref{Th:CanonBasis}
guarantees that for any $a=(a_1, \dots, a_\ell)$, $a' =
(a'_1,\dots,a'_\ell) \in \ZZ_{\ge 0}^\ell$ with  $a'_1 > a_1$ we have
 \begin{gather*}
(u^{a'},b_a)_c = 0.
 \end{gather*}
Equivalently, taking into account that $(e_2^{r+1}
u,b_a)_c=(u,L^{r+1}(b_a))_c$ for all $r\ge 0$, we obtain
$(S(V)^W,\nabla_{e_2}^{a_1+1}(b_a))_c = 0$. Since the form
\eqref{Eq:DefForm} is nondegenerate, we obtain
 \begin{gather*}
L^{a_1+1}(b_a) = 0
 \end{gather*}
for all $a\in \ZZ_{\ge 0}^\ell$. This proves that ${\bf B}_r = \{b_a
\mid a_1\le r\}$ is a (linearly independent) subset of the kernel of
$\left. L^{r+1}\right|_{S_c(V)^W}$. On the other hand, since $e_2$ and $L$
form a representation of $\mathfrak{sl}_2$ by Proposition \ref{Pp:SL2}, we
obtain isomorphisms of graded spaces:
 \begin{gather*}
S_c(V)^W \cong \CC(c)[e_2]\otimes {\mathcal K},\\
\name{Ker} \left.L^{n+1}\right|_{S_c(V)^W} = (\sum_{r=0}^n \CC(c) \cdot
e_2^r) \otimes {\mathcal K}
 \end{gather*}
where ${\mathcal K}$ is the kernel of $\left.L\right|_{S_c(V)^W}$. In
particular, the Hilbert series of ${\mathcal K}$ is
$\prod\limits_{k=2}^\ell \frac{1}{1-t^{d_k}}$, so that
 \begin{gather*}
\dim ({\mathcal K}\cap S_c^d(V)) = \lmod{\bf B}_0\cap S^d_c(V)\rmod,
 \end{gather*}
hence
 \begin{gather*}
\dim (\name{Ker} \left.L^{r+1}\right|_{S_c(V)^W} \cap S_c^d(V)) = \lmod{\bf
B}_r\cap S^d_c(V)\rmod
 \end{gather*}
for all $r\ge 0$. This, together with the inclusion ${\bf B}_r \subset
\name{Ker} \left.L^{r+1}\right|_{S_c(V)^W}$, proves that ${\bf B}_r$ is a
basis of $\name{Ker} \left.L^{n+1}\right|_{S_c(V)^W}$. Part
\ref{It:Kernel} is proved.

To prove \ref{It:MultE2}, denote
 \begin{gather*}
\tilde b_a = e_2^{a_1} b_{(0,a_2,\dots,a_\ell)}
 \end{gather*}
for each $a = (a_1,\dots,a_\ell)\in \ZZ_{\ge 0}^\ell$. Since $\tilde b_a$
satisf\/ies condition \ref{It:Lead} of Theorem \ref{Th:CanonBasis}, to
prove that $b_a = \tilde b_a$ it suf\/f\/ices to verify that the elements
$\tilde b_a$ satisfy condition \ref{It:Othog} of the same theorem. This
is equivalent to the elements $\tilde b_a$ being pairwise orthogonal, i.e.,
$(\tilde b_a,\tilde b_{a'})_c=0$ whenever $a\ne a'$. Part \ref{It:Kernel}
guarantees that both $b_{0,a_2,\dots,a_\ell}$ and
$b_{0,a'_2,\dots,a'_\ell}$ are in the kernel of $L$, so by Lemma~\ref{le:Orthog} we obtain:
 \begin{gather*}
(\tilde b_a,\tilde b_{a'})_c \in \delta_{a_1,a'_1} \cdot
(\tilde b_{(0,a_2,\dots,a_\ell)}, b_{(0,a'_2,\dots,a'_\ell)})_c \cdot
\CC(c) = \delta_{a,a'}\cdot \CC(c)
 \end{gather*}
because the elements $b_{(0,a_2,\dots,a_\ell)}$ and
$b_{(0,a'_2,\dots,a'_\ell)}$ of the canonical basis ${\bf B}$ are
orthogonal unless $(a_2,\dots,a_\ell)=(a'_2,\dots,a'_\ell)$. This proves
\ref{It:MultE2}.

Proposition \ref{pr:kernel of sl_2} is proved.
 \end{proof}

 \begin{proof}[Proof of Theorem \ref{Th:Elem}]
Take a homogeneous set $u_1, \dots, u_\ell$ of generators in $S(V)^W$, and
denote, as usual, $u^a \bydef u_1^{a_1} \cdots u_\ell^{a_\ell}$. By
def\/inition, $(u^a, e_{d_k}\c)_c = 0$ for all $a \in \Integer_{\ge 0}^\ell$
such that $a_k = \dots = a_\ell = 0$. For any such non-zero $a \in
\Integer_{\ge 0}^\ell$ let $s \le k-1$ be the largest index such that $a_s
\ne 0$ and let $a' \bydef a - \delta_s \in \Integer_{\ge 0}^\ell$ (we
denote $\delta_s = (0, \dots, 1, \dots,0)$, where  $1$ is in the $s$-th
position). Then $(u^{a'}, \nabla_{u_s} (e_{d_k}\c))_c = 0$. In particular,
the element $\nabla_{u_s} (e_{d_k}\c)$ is $(\cdot,\cdot)_c$-orthogonal to
all of the monomials $u^{a'}$ such that $\deg u^{a'}  = d_k - d_s$. Since
the form $(\cdot,\cdot)_c$ is nondegenerate for generic $c$, this implies
$\nabla_{u_s} (e_{d_k}\c) = 0$.

By def\/inition $e_{d_k}\c \equiv u_k \mod S_c(V)_{\prec \delta_k}$ (see~\eqref{eq:SVpreca} and the property~\ref{It:Lead} of the canonical
basis~$b_a$), which implies that $(e_{d_1}\c)^{a_1} \dots
(e_{d_\ell}\c)^{a_\ell} \equiv u^a \bmod S_c(V)_{\prec a}$. Indeed, if for
every $k = 1, \dots, h$
 \begin{gather*}
e_{d_k}\c = u_k + \sum \beta_{p_1,
\dots, p_{k-1}} u_1^{p_1} \cdots u_{k-1}^{p_{k-1}}
 \end{gather*}
for some $\beta_{p_1, \dots, p_{k-1}} \in \Complex(c)$ then
 \begin{gather*}
(e_{d_1}\c)^{a_1} \dots (e_{d_\ell}\c)^{a_\ell} = u_1^{a_1} \cdots
u_\ell^{a_\ell} + \sum_{p_1, \dots, p_{\ell-1}}
\sum_{p_\ell=1}^{a_{\ell}-1} \gamma_{p_1, \dots, p_\ell} u_1^{p_1} \cdots
u_\ell^{p_\ell}
 \end{gather*}
for some $\gamma_{p_1, \dots, p_\ell} \in \Complex(c)$. Here $(p_1,
\dots,p_\ell) \prec (a_1, \dots, a_\ell)$, because the last non-zero
coordinate of $(a_1-p_1, \dots, a_\ell-p_\ell)$ is positive. Thus, the
monomials $(e_{d_1}\c)^{a_1} \cdots (e_{d_\ell}\c)^{a_\ell}$ for all $a_1,
\dots, a_\ell \in \Integer_{\ge 0}$ form a basis of $S_c(V)^W$, and part
\ref{It:Generates} of Theorem \ref{Th:Elem} is proved.

To prove part \ref{It:EqForEDk} let $\mu_k\c \in S_c(V)^W$ be any
element satisfying its conditions. Then for any $a \prec \delta_k$ one has
$(u^a,\mu_k\c)_c = \nabla_{u_1}^{a_1} \cdots \nabla_{u_{k-1}}^{a_{k-1}}
(\mu_k\c) = 0$, hence $\mu_k\c = {\rm const} \cdot e_{d_k}\c + \sum_{\delta_k
\prec q} \beta_q b_q$ for some $\beta_q \in \Complex(c)$. The element
$b_q$ is homogeneous with $\deg b_q = \sum_{i=1}^\ell d_i q_i$. Since $d_1
< \dots < d_k < \dots < d_\ell$, one has $\deg b_q > d_k$ for every $q$.
Therefore, the homogeneity condition implies that $\mu_k\c = {\rm const} \cdot
e_{d_k}\c$. Part \ref{It:EqForEDk} is proved.

Part \ref{It:EqForKh}, given here for completeness, is proved in our
previous paper \cite{TAMS}.

Theorem \ref{Th:Elem} is proved.
 \end{proof}

\subsection{Proof of Theorem \ref{Th:Iwasaki}} \label{subsec:Iwasaki}

The argument follows almost literally the original proof for $c=0$ given in
\cite{Iwasaki}; we put it here mostly for reader's convenience.

Abbreviate $\nabla_{x_i}$ as $\nabla_i$; also take $x \bydef (x_1, \dots,
x_n)$ and $\nabla \bydef (\nabla_1, \dots, \nabla_n)$.

Let $c$ be generic.

 \begin{lemma}
A function $f$ can be represented as $(P(\nabla)) (\Delta(x))$ for some
polynomial $P$ if and only if $(e_k(\nabla))(f(x)) = 0$ for all $k$.
 \end{lemma}

 \begin{proof}
For $c = 0$ (when $\nabla_i = \partial_{x_i}$) it is a theorem due to
Steinberg \cite{Steinberg}. For $c$ generic the Dunkl operators are
conjugate to dif\/ferentiations by means of some intertwining operator $B$
(see~\cite{DunklJeuOpdam}): $\nabla_i = B^{-1} \partial_{x_i} B$, and
therefore $e_k(\nabla)(f(x)) = 0$ is equivalent to $e_k(\partial) B(f(x)) =
0$. By the Steinberg's theorem it means that $f = B^{-1} P(\partial) \Delta
= P(\nabla) B^{-1} \Delta$. But since $\Delta$ is skew-symmetric, $\nabla_i
\Delta$ is proportional to $\partial_{x_i} \Delta$, and therefore $B^{-1}
\Delta = \lambda \Delta$ for some constant $\lambda$. The lemma is proved.
 \end{proof}

Consider now the polynomials $f_i = \nabla_i (e_n\c)$. Obviously,
$(e_k(\nabla)) (f_i) = 0$ for all $k$, and therefore $f_i =
(g_i(\nabla))(\Delta(x))$. Without loss of generality, $g_i$ can be taken
skew-invariant with respect to the subgroup $G_i \subset S_n$ of
permutations leaving $i$ f\/ixed; from degree considerations one obtains
$g_i(x) = {\rm const} \cdot \Delta(x_1, \dots, \widehat{x_i}, \dots, x_n)$.
Euler's formula $e_n\c = n \cdot \sum_{i=1}^n x_i f_i$ implies now~\eqref{Eq:Iwa} for $k = n$.

Let $p_m(y_1, \dots, y_s) \bydef y_1^m + \dots + y_s^m$. Since
\eqref{Eq:Iwa} for $k = n$ is proved, one has
 \begin{gather*}
(p_i(\nabla_{j_1}, \dots, \nabla_{j_k})) (f_{j_1, \dots, j_k}(x)) = {\rm const}
\cdot \delta_{ik},
 \end{gather*}
where $1 \le j_1 < \dots < j_k \le n$ and $f_{j_1, \dots, j_k} \bydef
e_k\c(x_{j_1}, \dots, x_{j_k})$. Thus, if $F_k$ is the right-hand side of
\eqref{Eq:Iwa}, one has
 \begin{gather*}
p_i(\nabla) F_k = {\rm const} \cdot \sum_{1 \le j_1 < \dots < j_k \le n}
(p_i(\nabla))(f_{j_1, \dots, j_k}(x)) \\
\phantom{p_i(\nabla) F_k}{} = {\rm const} \cdot \sum_{1 \le j_1 < \dots < j_k \le n} (p_i(\nabla_{j_1},
\dots, \nabla_{j_k}))(f_{j_1, \dots, j_k}(x)) = {\rm const} \cdot \delta_{ik},
 \end{gather*}
Hence, $F_k = e_k\c$.

Theorem \ref{Th:Iwasaki} is proved.

\subsection{Proof of Theorem \ref{Th:Main}}\label{SSec:PrMain}
We need the following result.

 \begin{proposition}\label{pr:well defined at 1/n}
The canonical basis ${\bf B}$ is well-defined at $h_c=1$ $($e.g., at $c=1/h)$
and ${\bf B}\setminus\{1\}$ is orthogonal with respect to the form
$\overline \Phi_c$ $($see Theorem~{\rm \ref{Th:NonDeg}}$)$.
 \end{proposition}

 \begin{proof}
Indeed, in the def\/inition of $b_a$ in Theorem \ref{Th:CanonBasis} we can
take  $u_1, \dots, u_\ell$ to be a homogeneous generating set of
$S(V_\Real)^W_+$. Let us prove that for each $a \in \Integer_{\ge 0}^\ell
\setminus \{0\}$ the coef\/f\/icients $c_{a,a'}$ of the expansion $b_a =
\sum_{a' \preceq a} c_{a,a'} u^a$ have no poles at $h_c=1$. Suppose the
contrary. Then there exists an exponent $\ell > 0$ such that $\tilde b_a =
(1-h_c)^\ell b_a$ is regular at $h_c = 1$ and $\pi(\tilde b_a) \ne 0$,
where $\pi:{\mathcal A} \to {\mathcal A}/(1-h_c)$ is the canonical
homomorphism (${\mathcal A}$ is as in Proposition \ref{pr:non-isotropic
Phi}). Since $\ell > 0$, we see that $\pi(\tilde b_a)\in U \bydef
\sum_{a'\prec a,a'\ne 0} \kk \cdot u^{a'}$, where  $\kk = {\mathcal
A}/(1-h_c)$. But since $\Phi_c(b_a,u^{a'}) = 0$ for all $a'\prec a$, we see
that $\overline \Phi_c(\pi(\tilde b_a),U) = 0$ which contradicts Theorem
\ref{Th:NonDeg}\ref{It:PhiBarSub}. The contradiction obtained proves that
$\ell=0$, i.e., $\pi(b_a)$ is well-def\/ined. The orthogonality of these
elements is obvious.
 \end{proof}

Let $W = S_n$, $V = \Complex^n$. Here $d_k = k$, $k = 1, \dots, n$.  Recall
that $e_k = e_k(x_1, \dots, x_n) \in \Complex[x_1,\dots,x_n]$ is the $k$-th
elementary symmetric polynomial and $\bare_k = e_k(x_1 - \frac{e_1(x)}{n},
\dots, x_n - \frac{e_1(x)}{n})$. The elements $\bare_k$, $k = 2,\dots, n$
generate a subalgebra of $S(V)^{S_n}$ isomorphic to $S(V')^{S_n}$, where
$V' = \{x \in V \mid \sum_i x_i=0\}$.

 \begin{lemma} \label{Lm:NablaPE}
Using the notation  $p_r \bydef y_1^r + \dots + y_n^r$, we obtain
 \begin{gather*}
\nabla_{p_r} (\bare_k) = (-1)^r (n-k+1)\cdots(n-k+r) \frac{(1-nc)^r -
(1-nc)}{n^r c} \bare_{k-r}.
 \end{gather*}
 \end{lemma}

 \begin{proof}
Denote
 \begin{gather*}
Q_{k,i}(x) \bydef e_k\biggl(x_1 - \frac{e_1(x)}{n}, \dots, \widehat{x_i -
\frac{e_1(x)}{n}}, \dots, x_n - \frac{e_1(x)}{n}\biggr)
 \end{gather*}
(the $i$-th argument omitted). Easy calculations show that
 \begin{gather}
\sum_{i = 1}^n Q_{k,i} = (n-k) \bare_k, \label{Eq:SumP}\\
\pder{Q_{k,i}}{x_i} = -\frac{n-k}{n} Q_{k-1,i},\nonumber\\
\pder{\bare_k}{x_i} = \nabla_{y_i} (\bare_k) = Q_{k-1,i} -
\frac{n-k+1}{n} \bare_{k-1}.\label{Eq:NablaE}
 \end{gather}
Clearly,
 \begin{gather*}
\sum_{j \ne i} \frac{(ij) Q_{k,i} - Q_{k,i}}{x_i - x_j} = (n-k)
Q_{k-1,i},
 \end{gather*}
and therefore the Dunkl operator
 \begin{gather*}
\nabla_{y_i} (Q_{k,i}) = (n-k)\left(c - \frac{1}{n}\right) Q_{k-1,i}.
 \end{gather*}
An induction on $r$ (where \eqref{Eq:NablaE} is the base) gives then
 \begin{gather*}
\nabla_{y_i}^r (\bare_k) = (n-k+1) \cdots  (n-k+r-1) \frac{(nc-1)^r -
(-1)^r}{n^r c} Q_{k-r,i} \\
\phantom{\nabla_{y_i}^r (\bare_k) =}{} + (-1)^r \frac{(n-k+1) \cdots (n-k+r)}{n^r} \bare_{k-r}.
 \end{gather*}
Summation over $i$ using \eqref{Eq:SumP} f\/inishes the proof.
 \end{proof}

Lemma \ref{Lm:NablaPE} implies that $\nabla_{p_r}\c (\bare_k) = 0$ for $c =
1/n$ and all $r <k$. Since $p_1, \dots, p_n$ generate $\Complex[x_1, \dots,
x_n]^{S_n}$, we obtain
 \begin{gather*}
\nabla_P\c (\overline e_k) = 0
 \end{gather*}
for $c = 1/n$ and any homogeneous symmetric polynomial $P$ of
degree $\deg P < k$.

On the other hand, by Theorem \ref{Th:Elem}\ref{It:EqForEDk} one has
$\nabla_P (e_k\c) = 0$ for any homogeneous symmetric polynomial $P$ of
degree $\deg P < k$. Therefore,
 \begin{gather*}
\overline\Phi_{1/n}(\bare^{a'}, \bare_k) = 0 =
\overline\Phi_{1/n}\big(\bare^{a'}, e_k^{(1/n)}\big)
 \end{gather*}
for any $a'\prec \delta_k$. Here we abbreviated $\overline\Phi_{1/n} \bydef
\overline \Phi_c$ with $c = 1/n$ in the notation of Theorem \ref{Th:NonDeg}
and $e_k^{(1/n)} \bydef \lim_{c\to 1/h} e_k\c$, which is well-def\/ined by
Proposition \ref{pr:well defined at 1/n}. Consequently,
 \begin{gather*}
\overline\Phi_{1/n}(U, e_k^{(1/n)}-\alpha \bare_k) = 0
 \end{gather*}
for all $a'\prec \delta_k$, $\alpha\in \Complex^\times$, where $U =
\sum_{a'\prec \delta_k} \Complex \cdot {\overline e}^{\,a'}$. Therefore,
$e_k^{(1/n)} = \alpha \bare_k$ because, on the one hand,
$e_k^{(1/n)}-\alpha \bare_k\in U$ for some $\alpha \ne 0$, and on the other
hand, the restriction of  $\overline\Phi_{1/n}$ to $U$ is non-degenerate by
Theorem \ref{Th:NonDeg}\ref{It:PhiBarSub}.

Theorem \ref{Th:Main} is proved.

\subsection{Canonical invariants of dihedral groups and proof of Theorem
\ref{Th:DihCconst}}\label{sect:dihedral}

Throughout the section we deal with the dihedral group
 \begin{gather*}
W=I_2(m)=\langle s_0^2 = s_1^2 = (s_0 s_1)^m = 1\rangle
 \end{gather*}
of order $2m$.

We denote by $\{z,\overline z\}$ the basis of $V$ such that
 \begin{gather}\label{Eq:Action}
s_j(z) = -\zeta^j \bar z, \qquad s_j(\bar z) = -\zeta^{-j} z.
 \end{gather}
for $j = 0,1$, where $\zeta = e^{2\pi i/m}$ is an $m$-th primitive root of
unity.

We also denote by $e_2$, $e_m$ the generators of $S(V)^W$ given by
 \begin{gather*}
e_2 = z\overline z, \qquad e_m = z^m + \overline z^m.
 \end{gather*}

 \begin{lemma}\label{le:Dunkl Laplacian}
The restriction of the Dunkl Laplacian $L$ to $S(V)^W=\Complex[e_2,e_m]$
for $W = I_2(m)$ equals:
 \begin{gather*}
L = e_2\partial^2_{e_2} + m e_m\partial_{e_2}\partial_{e_m} +
m^2e_2^{m-1}\partial^2_{e_m} + \left(1-\frac{m}{2}C\right)\partial_{e_2} +
\frac{m^2}{2} \delta e_2^{m/2-1} \partial_{e_m},
 \end{gather*}
where $C \bydef c(s_1)+ c(s_2)$ and $\delta \bydef c(s_2)-c(s_1)$ $($so that
$\delta = 0$ when $m$ is odd$)$; $\partial_{e_2}$ and $\partial_{e_m}$ mean
here differentiation with respect to $e_2$ and $e_m$, respectively, in the
ring $\Complex[e_2, e_m]$.
 \end{lemma}

 \begin{proof}
Clearly, $s \circ \partial_{y} = \partial_{s^*(y)}\circ s$ for any linear
automorphism $s$ of $V$ and any $y\in V^*$, where $\partial_{y}:S(V)\to
S(V)$ is the directional derivative. So, $L$ can be rewritten in the form
$L = \sum_i D_i s_i$ where $D_i$ are dif\/ferential operators (of order at
most $2$) and $s_i$ are ref\/lections. Thus, on the space of {\em invariant}
functions $L$ is a second order dif\/ferential operator. To determine its
coef\/f\/icients it suf\/f\/ices to compute $L$ on monomials of degree $1$ and $2$
in $e_2$, $e_m$. This is done in \cite{DuDih}; see also~\cite{TAMS}.
 \end{proof}

 \begin{corollary}\label{Cr:ViaD}
In the notation of Lemma {\rm \ref{le:Dunkl Laplacian}}, we have
 \begin{gather*}
\frac{4}{m^2}e_2L = {\mathcal E}^2 - C{\mathcal E} +
\left(4\frac{e_2^m}{e_m^2}-1\right)({\mathcal D}^2 - {\mathcal D}) + \left(C-1 + 2 \delta
\frac{e_2^{m/2}}{e_m}\right) {\mathcal D}
 \end{gather*}
where ${\mathcal E} \bydef \frac{2}{m}e_2\partial_{e_2}+e_m\partial_{e_m}$
is a multiple of the Euler derivation, and ${\mathcal D} \bydef
e_m\partial_{e_m}$.
 \end{corollary}

 \begin{proof}
Since $\partial_{e_2} e_2 = e_2\partial_{e_2} + 1$, one obtains
 \begin{gather*}
4e_2L = (2e_2\partial_{e_2})^2 + 4m
(e_2\partial_{e_2})(e_m\partial_{e_m}) + 4m^2e_2^m\partial^2_{e_m} -
2mCe_2\partial_{e_2} + 2m^2 \delta e_2^{m/2} \partial_{e_m}.
 \end{gather*}
Equalities $(2e_2\partial_{e_2})^2 + 4m
(e_2\partial_{e_2})(e_m\partial_{e_m}) = m^2({\mathcal E}^2-{\mathcal
D}^2)$ and $\partial^2_{e_m} = \frac{1}{e_m^2}({\mathcal D}^2 - {\mathcal
D})$ imply now that
 \begin{gather*}
\frac{4}{m^2} e_2L = {\mathcal E}^2 - {\mathcal D}^2 +
4\frac{e_2^m}{e_m^2} ({\mathcal D}^2-{\mathcal D}) - C({\mathcal
E}-{\mathcal D}) + 2\delta\frac{e_2^{m/2}}{e_m} {\mathcal D},
 \end{gather*}
and the corollary follows.
 \end{proof}

\begin{corollary}[of Corollary \ref{Cr:ViaD}] \label{Cr:Subst}
Let $f = f(x,u) \in \Complex[[x,u]]$. Then
 \begin{gather*}
\frac{4}{m^2} e_2 L \bigg(f\bigg(\frac{e_m}{e_2^{m/2}}, e_2^{m/2}t\bigg)\bigg) =
L_{x,u}(f(x,u))|_{x = \frac{e_m}{e_2^{m/2}},u=e_2^{m/2}t}
 \end{gather*}
where
 \begin{gather*}
L_{x,u} = u^2\partial^2_u - (C-1)u \partial_u + (4-x^2)\partial_x^2 +
((C-1)x + 2\delta)\partial_x.
 \end{gather*}
\end{corollary}

\begin{proof} Follows from
  \begin{gather*}
{\mathcal D}\biggl(f\biggl(\frac{e_m}{e_2^{m/2}}, e_2^{m/2}t\biggr)\biggr) =
\frac{e_m}{e_2^{m/2}}f_x\bigl(e_m/e_2^{m/2}, e_2^{m/2}t\bigr) =
(x\partial_x)\left.f(x,u)\right|_{x =
\frac{e_{m\vphantom{a_b}}}{e_2^{m/2}},\,u=e_2^{m/2}t}
 \end{gather*}
and
  \begin{gather*}
{\mathcal E}\biggl(\biggl(\frac{e_m}{e_2^{m/2}}\biggr)^p \bigl(e_2^{m/2}
t\bigr)^q\biggr)  = \biggl(\frac{2}{m}e_2\partial_{e_2} +
e_m\partial_{e_m}\biggr) \big(e_m^p e_2^{(q-p)m/2} t^q\big)\\
\phantom{{\mathcal E}\biggl(\bigl(\frac{e_m}{e_2^{m/2}}\bigr)^p \bigl(e_2^{m/2}
t\bigr)^q\biggr) }{}
= q\big(e_m^p e_2^{(q-p)m/2} t^q\big) = (u\partial_u) \left.(x^p u^q)\right|_{x =
\frac{e_{m\vphantom{a_b}}}{e_2^{m/2}},\,u=e_2^{m/2}t}.\tag*{\qed}
 \end{gather*}\renewcommand{\qed}{}
 \end{proof}

 \begin{proof}[Proof of Theorem \ref{Th:DihCconst}]
In view of Corollary \ref{Cr:Subst}, for $c = {\rm const}$ it suf\/f\/ices to prove
that $L_{x,u}(p^c) = 0$, where $p = 1+xu+u^2$. Indeed,
 \begin{gather*}
u\partial_u(p^c) = cu(x+2u)p^{c-1}=2cp^c-c(2+ux)p^{c-1},\\
((u\partial_u)^2 - 2cu\partial_u)(p^c) =
u\partial_u(u\partial_u-2c)(p^c) = -cu\partial_u((2+ux)p^{c-1})\\
\phantom{((u\partial_u)^2 - 2cu\partial_u)(p^c)}{} = -cuxp^{c-1} - c(c-1)u(2+ux)(x+2u)p^{c-2},\\
\partial_x(p^c) = cup^{c-1},\\
\partial_x^2(p^c) = c(c-1)u^2p^{c-2}.
 \end{gather*}
For $c={\rm const}$ one has $C = 2c$ and $\delta=0$, so that
 \begin{gather*}
c^{-1} p^{2-c} L_{x,u}(p^c) = -uxp - (c-1)u(2+ux)(x+2u) + (4-x^2)(c-1)u^2
+ (2c-1)xup\\
\phantom{c^{-1} p^{2-c} L_{x,u}(p^c)}{} = (c-1)u(2xp-(2+ux)(x+2u) + (4-x^2)u)\\
\phantom{c^{-1} p^{2-c} L_{x,u}(p^c)}{}= 0.
 \end{gather*}
This proves part \ref{It:Cconst} of Theorem \ref{Th:DihCconst}.

To prove part \ref{It:CNotConst} def\/ine $I_r^{(a,b)}(y,u)\in \CC[[y,u]]$
for each $r$, $a$, $b$ by
 \begin{gather*}
I_r^{(a,b)} = \int_0^1 s^{a+b}(1-s)^{-b-1-r}
\left(1-s+us\left(1-\frac{s}{2}(1-y)\right)\right)^r ds.
 \end{gather*}
Clearly, the right-hand side of \eqref{eq:integral generating series m
even} equals $I_r^{(a,b)}(y,u)$ with $a = -(C+\delta+1)/2$, $b =
-(C-\delta+1)/2$, $r = -b-1$, $u = e_2^{m/2}t$, and $y =
\frac{e_m}{2e_2^{m/2}}$ (so that $\frac{1-y}{2} =
-\frac{e'_m}{e_2^{m/2}}$).

Therefore, in view of Corollary \ref{Cr:Subst}, it suf\/f\/ices to prove that
 \begin{gather}\label{eq:Lxu kills Irxu}
L_{x,u}(I_r^{(a,b)}(x/2,u)) = 0
 \end{gather}
for all $r$, where $a$, $b$ are as above, and to determine the normalizing
coef\/f\/icients $n_k(c)$.

Recall that the $n$-th Jacobi polynomial $P_n^{(a,b)}(y)$ is given by
 \begin{gather}\label{eq:Jacobi polynomial}
P_n^{(a,b)}(y) = \frac{\Gamma(a+n+1)}{n!\Gamma(n+a+b+1)\Gamma(-n-b)}
\int_0^1 \! s^{n+a+b} (1-s)^{-n-b-1} \left(1-\frac{s}{2}(1-y)\right)^n ds\!\!
 \end{gather}
(with the analytic continuation to all $a,b\in \CC$). Thus one has
 \begin{gather*}
I_r^{(a,b)}(y,u) = \int_0^1 s^{a+b} (1-s)^{-b-1} \sum_{n=0}^\infty
\frac{\Gamma(r+1)}{n!\Gamma(r-n+1)} u^n s^n (1-s)^{-n} \left(1 -
\frac{s}{2}(1-y)\right)^n ds \\
\phantom{I_r^{(a,b)}(y,u)}{} = \sum_{n=0}^\infty \frac{\Gamma(r+1)}{n!\Gamma(r-n+1)} u^n \int_0^1
s^{a+b+n} (1-s)^{-b-1-n} \left(1 - \frac{s}{2}(1-y)\right)^n ds\\
\phantom{I_r^{(a,b)}(y,u)}{}= \sum_{n=0}^\infty \frac{\Gamma(r+1)\Gamma(n+a+b+1)\Gamma(-n-b)}
{\Gamma(r-n+1)\Gamma(a+n+1)} P_n^{(a,b)}(y) u^n\\
\phantom{I_r^{(a,b)}(y,u)}{}\bydef \sum_{n=0}^\infty q_n(r,a,b) P_n^{(a,b)}(y) u^n.
 \end{gather*}

The Jacobi polynomial $P_n^{(a,b)}(y)$ belongs to the kernel of the
dif\/ferential operator
 \begin{gather*}
J^{(a,b)} \bydef (1-y^2)\partial_y^2 + (b-a-(a+b+2)y)\partial_y +
n(n+a+b+1)
 \end{gather*}
(see e.g.\ \cite{SpecFunc} for proof). Therefore, $I^{(a,b)}(y,u) =
\sum_{n=0}^\infty q_n(r,a,b) P_n^{(a,b)}(y) u^n$ satisf\/ies
 \begin{gather}\label{eq:vanishing of integrals}
\tilde L_{y,u}^{(a,b)}(I_r^{(a,b)}(y,u))=0,
 \end{gather}
where
$\tilde L_{y,u}^{(a,b)} = (1-y^2)\partial_y^2 + (b-a-(a+b+2)y)\partial_y +
u^2\partial_u^2 + (a+b+2)u\partial_u$.

Take now $a=-(C+\delta+1)/2$, $b=-(C-\delta+1)/2$, so that $a+b+2=-(C-1)$,
$b-a=\delta$. One has then
 \begin{gather*}
\tilde L_{y,u}^{(a,b)} = (1-y^2)\partial_y^2 + (\delta+(C-1)y)\partial_y +
u^2\partial_u^2 - (C-1)u\partial_u\\
\phantom{\tilde L_{y,u}^{(a,b)}}{}= (4-(2y)^2)\partial_{2y}^2 + (2\delta+(C-1)(2y))\partial_{2y} +
u^2\partial_u^2 - (C-1)u\partial_u\\
\phantom{\tilde L_{y,u}^{(a,b)}}{}= L_{2y,u}.
 \end{gather*}
Therefore, \eqref{eq:vanishing of integrals} implies \eqref{eq:Lxu kills
Irxu}.

To f\/inish the proof it remains to f\/ind the value of the normalization
coef\/f\/icients $n_k=n_k(c)$ in \eqref{eq:integral generating series m even}.
To do this, substitute $e_2 = 0$ into \eqref{eq:integral generating series
m even}, so that $e_m' = e_m/4$. Under this specialization, $b_{(0,k)}$
becomes a (complex) multiple of $e_m^k$ and the right-hand side of
\eqref{eq:integral generating series m even} becomes (with the abbreviation
$\alpha=\frac{C-\delta-1}{2}$):
 \begin{gather*}
\int_0^1 (1-s)^\alpha s^{-C-1} (1 + \frac{ts^2}{4(1-s)} e_m)^\alpha ds \\
\qquad{} = \sum_{k=0}^\infty \binom{\alpha}{k}t^k\left(\frac{e_m}{4}\right)^k
\int_0^1 (1-s)^{\alpha-k} s^{-C+2k-1} ds\\
\qquad{} = \sum_{k=0}^\infty  \binom{\alpha}{k} \frac{\Gamma(\alpha-k+1)
\Gamma(-C+2k)} {\Gamma(\alpha+k-C+1)} t^k \left(\frac{e_m}{4}\right)^k\\
\qquad{}= \sum_{k=0}^\infty \frac{\Gamma(\alpha+1) \Gamma(-C+2k)} {k!\cdot
\Gamma(\alpha+k-C+1)} t^k \left(\frac{e_m}{4}\right)^k\\
\qquad{} =\sum_{k=0}^\infty \frac{\Gamma(\frac{C-\delta+1}{2}) \Gamma(2k-C)}
{\Gamma(k - \frac{C+\delta-1}{2})} t^k \frac{e_m^k}{4^k k!}.
 \end{gather*}
This proves \eqref{eq:integral generating series m even} and f\/inishes the
proof of part \ref{It:CNotConst} of Theorem \ref{Th:DihCconst}.
 \end{proof}

\subsection*{Acknowledgments}

The work
was partially supported by the CRDF grant RUM1-2895-MO-07. The second
author was also supported by the INTAS grant 05-7805, RFBR grants
08-01-00110-a and NSh-709.2008.1, and the HSE Scientif\/ic Foundation grant
08-01-0019.

The authors are grateful to P.~Etingof, M.~Feigin, A.~Samokhin, and
Y.~Xu for valuable discussions. The second author wishes to thank the
University of Oregon, where most of this work was carried out, for its warm
hospitality.

\pdfbookmark[1]{References}{ref}
\LastPageEnding

\end{document}